\documentclass[a4paper,leqno,10pt]{article}
\usepackage{amssymb}
\usepackage{amsmath}
\usepackage[german]{babel}
\begin{document}

\def\ord{\mathop{\rm ord}}

\def\refname{ }
\begin{center}
\begin{Large}
\textbf{Pad\'e approximations of a class of  G-functions and some applications}
\end{Large}
\vskip0.4cm
\sc{Keijo V\"a\"an\"anen}
\end{center}
\vskip0.6cm

\begin{center}
Abstract
\end{center}

We construct explicitly Pad\'{e} approximations of the second kind for a special class of $G$-functions. These are then applied to prove a Baker-type lower bound for linear forms in the $p$-adic values of these functions. Moreover, we consider restricted rational approximations of the values of these functions in the real case.\\

\noindent
2010 Mathematics Subject Classification: 11J13, 11J61, 11J82

\noindent
Keywords: $G$-function, Pad\'e approximation, $p$-adic linear form, Baker-type lower bound

\section{Introduction}

In the present paper we shall consider the $G$-functions
\begin{equation}\label{e1}
\varphi_j(z) = \sum_{n=0}^\infty\frac{(\alpha_j)_n}{(\alpha_j+\alpha_0)_n}z^n, \quad j = 1,\ldots, m,
\end{equation}
where $\alpha_0, \alpha_1,\ldots, \alpha_m$ are positive rational numbers satisfying $\alpha_i - \alpha_j \notin \mathbb{Z}, \ 1 \leq i < j \leq m$, and $(\alpha)_0 = 1, (\alpha)_n = \alpha(\alpha+1)\cdots(\alpha+n-1), n \geq 1$. We first construct Pad\'{e} approximation polynomials of the second kind for the series (\ref{e1}). Then we give two applications, in the first one Baker-type lower bounds are obtained for the $p$-adic values $\varphi_j(a/b)$, where rational $a/b \neq 0$ has small $p$-adic value. The second application studies analogously to \cite{FR} special restricted rational approximations of the real values $\varphi_j(a/b)$ with some rational $a/b \neq 0$.

A basic problem in Diophantine approximations is to find lower bounds for the absolute values of linear forms $\ell_1\theta_1+\cdots+\ell_m\theta_m, (\ell_1,\ldots,\ell_m)\in\mathbb{Z}^m\setminus\{0\}$, in given linearly independent numbers $\theta_1,\ldots,\theta_m$. Often such bounds, linear independence measures, are given in terms of $h=\max \{\left|\ell_1\right|,\ldots,\left|\ell_m\right|\}$, but also more refined bounds in terms of each individual $h_j=\max \{1,\left|\ell_j\right|\}$ are of interest. These refined bounds are called Baker-type bounds, since Baker \cite{B} was the first to obtain such a bound in the case $\theta_j=e^{\alpha_j}$ with distinct rational $\alpha_1,\ldots,\alpha_m$. More precisely, he proved that 
\[
\left|\ell_1e^{\alpha_1} + \cdots + \ell_me^{\alpha_m}\right| > h^{1 - c_0/\sqrt{\log\log h}}\prod_{j=1}^m h_j^{-1},
\]
for all $h \geq c_1 > e$, where $c_0, c_1$ are positive constants depending on $\alpha_j$. These constants were made completely explicit in Mahler \cite{M1}, for further improvements see also \cite{S} and \cite{ELM}. Baker's proof used essentially Siegel's method with a new idea in the construction of an auxiliary function, a Pad\'{e} type approximation of the first kind for the functions $e^{\alpha_jz}$, obtained by using Siegel's lemma. This is a special linear form in $e^{\alpha_jz}$ with polynomial coefficients having a zero of high order at $z=0$. After that the same idea was used to study other $E$- and $G$-functions satisfying linear differential equations of first order with rational coefficients, see for example \cite{M} and \cite{V}. Then, in an important and deep paper \cite{Z}, Zudilin was able to obtain a similar result for the values of a class of E-functions satisfying a system of homogeneous linear differential equations with rational coefficients, in this general result the term $\sqrt{\log\log h}$ in the bound is replaced by $(\log\log h)^{1/(m^2-m+2)}$. In all these papers the approximations are obtained by applying Siegel's lemma.

The first explicit Pad\'e approximations of the first kind suitable for Baker-type bounds were given by Fel'dman in \cite{F} and \cite{F1}, in particular \cite{F} contains a result for some $E$-functions, where $\sqrt{\log\log}$ in the estimate is improved to $\log\log$. This result is slightly improved in \cite{V1} by using explicit Pad\'e approximations of the second kind, in other words, simultaneous rational approximations to the functions under consideration. In \cite{F1} Fel'dman considered linear forms in the values of the functions $\varphi_j(z)$ in the special case $\alpha_0 = 1$, and then Sorokin improved and generalized Fel'dman's result by proving the following theorem \cite[Theorem 2]{So}.\\

\noindent
\textbf{Theorem} (Sorokin). \textit{Let $K$ denote $\mathbb{Q}$ or an imaginary quadratic field, and assume that $\beta = a/b \neq 0$, where $a,b \in\mathbb{Z}_K$, the ring of integers of $K$. Further, let $\epsilon, 0 < \epsilon < 1$, be given. Then there exists a positive constant $c$ depending on $\alpha_0,\alpha_1,\ldots,\alpha_m$ and $\epsilon$ such that if
\begin{equation}\label{2}
\left|b\right|^\epsilon > c\left|a\right|^{2(m+1+m\epsilon)},
\end{equation}
then
\[
\left|\beta_0 + \beta_1\varphi_1(\beta) + \cdots + \beta_m\varphi_m(\beta)\right| > H^{-1-\epsilon}
\]
for all $(\beta_0,\beta_1,\ldots,\beta_m) \in \mathbb{Z}_K^{m+1}\setminus{\{\underline{0}\}}$ with $H = \prod_{j=1}^m h_j \geq H_0, h_j = \max\{1,\left|\beta_j\right|\}, \ j=1,\ldots,m$, where $H_0$ is a positive constant depending on $\alpha_0,\alpha_1,\ldots,\alpha_m, \beta$ and $\epsilon$.}\\

Sorokin's work is the first one, where appropriate Pad\'e approximations of the second kind are applied in this connection. Then in \cite{VZ} such type of construction obtained by Siegel's lemma was used to study certain $q$-series, for a refinement see also \cite{L}.

To prove the $p$-adic analogue of Sorokin's theorem we construct suitable Pad\'{e} approximations for the series $\varphi_j(z)$. For this, let $n_1,\ldots, n_m, N_1,\ldots, N_m$ be positive integers and $N = n_1+\cdots+n_m$. The polynomials $Q(z) \neq 0, P_1(z),\ldots, P_m(z)$ are called type $(N;N_1,\ldots,N_m)$ Pad\'{e} approximation polynomials of the second kind for the series $\varphi_j(z)$, if
\[
\deg Q(z) \leq N, \ \deg P_j(z) \leq N_j, \ \ord(Q(z)\varphi_j(z) - P_j(z)) \geq N_j+n_j+1, \quad j = 1,\ldots, m,
\]
where $\ord$ means the order of zero at $z = 0$.\\

\noindent
\textbf{Theorem 1}. \textit{If $N_j \geq N-1, j = 1,\ldots,m$, then the polynomials
\[
Q(z) = \sum_{k=0}^N a_kz^k,
\]
where
\begin{align}\label{e3}
a_N = 1, a_{N-k-1} =& \sum_{\ell=k}^{N-1}(-1)^{\ell+1}\frac{(\alpha_0-1)_{\ell-k}}{(\ell-k)!}\frac{(\alpha_0+\ell+1)_{N-\ell-1}}{(N-\ell-1)!}\prod_{j=1}^m\frac{(\alpha_j+\alpha_0+N_j-N+\ell+1)_{n_j}}{(\alpha_j+N_j-N+1)_{n_j}}, \\ &k = 0,1,\ldots,N-1.\notag
\end{align}
and
\[
P_j(z) = \sum_{\mu=0}^{N_j} c_{j\mu}z^\mu, \quad c_{j\mu} = \sum_{k=0}^{\min\left\{N,\mu\right\}}a_k\frac{(\alpha_j)_{\mu-k}}{(\alpha_j+\alpha_0)_{\mu-k}}, \quad \mu = 0,1,\ldots,N_j ; j=1,\ldots,m 
\]
are type $(N;N_1,\ldots,N_m)$ Pad\'{e} approximation polynomials of the second kind for the series $\varphi_j(z)$}.\\

We shall prove Theorem 1 by using a simple method given in Lemma 1 below. Note that the choice $N_1 = \cdots = N_m = N-1$ corresponds the approximations of \cite[Proposition 2]{So}, but for our applications it is essential to be able to use non-diagonal approximations with distinct $N_j$.

For our arithmetical results we introduce some notations. Let $\alpha_j =r_j/s_j, \ j = 0,1,\ldots,m$, and $\alpha_j+\alpha_0 = u_j/v_j, \ j = 1,\ldots,m$, where $r_j, s_j, u_j, v_j$ are positive integers satisfying $(r_j,s_j) = (u_j,v_j) = 1$, and denote $R = \max\left\{r_1,\ldots,r_m\right\}, S = \max\left\{s_1,\ldots,s_m\right\}, U = \max\left\{u_1,\ldots,u_m\right\}$ and $V = \max\left\{v_1,\ldots,v_m\right\}$. Further, let positive integers $d_j$ be defined by $s_0s_j = d_jv_j$, and let $s,v$ and $d$ be the least common multiples of $s_j, v_j$ and $d_j, \ j = 1,\ldots,m,$ respectively. Finally, we use the notation $\tilde{d} = d/(d,s_0)$, and put
\[
\epsilon(n) = \prod_{p\mid n} p^{1/(p-1)}
\]
for all integers $n \geq 1$. 

The series (\ref{e1}) converge $p$-adically in $\mathbb{Q}_p$, the $p$-adic completion of $\mathbb{Q}$, if $z$ is rational and satisfies $\left|z\right|_p < 2^{-\delta(2,p)}\left|s\right|_p$, where $\delta(2,p) = 1$, if $p=2$ and $s$ is even, and $\delta(2,p) = 0$ otherwise. As usual, the valuation is normalized by $\left|p\right|_p = p^{-1}$. So we may consider linear forms of $p$-adic numbers $\varphi_j(\beta)$, where $\beta = a/b \in \mathbb{Q}$ satisfies $0 < \left|\beta\right|_p < 2^{-\delta(2,p)}\left|s\right|_p$, note that this is simply the condition $0 < \left|\beta\right|_p < 1$, if $p\nmid s$. The following theorem gives a $p$-adic analogue of Sorokin's theorem and improves results obtained as a special case of \cite[Corollary 2]{V}.\\

\noindent
\textbf{Theorem 2}. \textit{Let $\alpha_0,\alpha_1,\ldots,\alpha_m$ satisfy the assumptions given after (\ref{e1}), and let $\epsilon, 0 < \epsilon < 1$, be given. Then there exists a positive constant $\tilde{c}$ depending on $\alpha_0,\alpha_1,\ldots,\alpha_m$ and $\epsilon$ such that if $\beta = a/b \in \mathbb{Q}\setminus\{0\}$ satisfies
\begin{equation}\label{e4}
(a,b) = 1, \left|a\right|^\epsilon > \tilde{c}\left|b\right|^{2(m+1+\epsilon)}, \left|a\right|_p < \min\left\{2^{-\delta(2,p)}\left|s\right|_p, \left|a\right|^{-1+\epsilon/(8(m+1))}\right\}, 
\end{equation}
then
\[
\left|\ell_0 + \ell_1\varphi_1(\beta) + \cdots + \ell_m\varphi_m(\beta)\right|_p > \tilde{H}^{-1-\epsilon}
\]
for all $(\ell_0,\ell_1,\ldots,\ell_m) \in \mathbb{Z}^{m+1}\setminus{\{\underline{0}\}}$ with $\tilde{H} = \prod_{j=0}^m h_j \geq \tilde{H}_0, h_0 = \max \{\left|\ell_0\right|,\left|\ell_1\right|,\ldots,\left|\ell_m\right|\}, h_j = \max\{1,\left|\ell_j\right|\}, \ j=1,\ldots,m$, where $\tilde{H}_0$ is a positive constant depending on $\alpha_0,\alpha_1,\ldots,\alpha_m, \beta$ and $\epsilon$.}\\

The explicit values of the constants $\tilde{c}$ and $\tilde{H}_0$ above are obtained with the choice $\tau = \epsilon/(m+1)$ from (\ref{16}) and (\ref{17}), respectively.

The last condition in (\ref{e4}) means that $a$ must be nearly a power of $p$. Our considerations below give some information also in the case of more general integers $a, \left|a\right| > 1$. Clearly the function values $\varphi_j(a)$ are defined in $\mathbb{Q}_p$ for all $p \mid a$, if $(a,s) = 1$. A linear relation
\[
\ell_0 + \ell_1\varphi_1(a) + \cdots + \ell_m\varphi_m(a) = 0
\]
with $(\ell_0,\ell_1,\ldots,\ell_m) \in \mathbb{Z}^{m+1}\setminus{\{\underline{0}\}}$ is called $a$-global if it holds for all $p \mid a$. It is known that such relations are not possible for sufficiently large $\left|a\right|$, see \cite{Bo}, and the following theorem gives an explicit bound for $\left|a\right|$.\\

\noindent
\textbf{Theorem 3}. \textit{Assume that $a$ is an integer satisfying $\left|a\right| > 1$ and $(a,s) = 1$. Let $\alpha_0,\alpha_1,\ldots,\alpha_m $ be as in Theorem 2. Then there are no $a$-global relations, if $\left|a\right| > C$, where}
\[
\log C = mS + (m+1)(3+\log(d\tilde{d}s_0\epsilon(s_0)\epsilon(s)^2\epsilon(v)) + 2s_0 + (m+1)V).
\]

As a second application of Theorem 1 we consider approximations of real values of $\varphi_j(a/b)$ at rational points $a/b \neq 0$ by rational numbers with denominators of the form $Bb^M$, where $B, M$ are positive integers and $B$ is fixed. Such studies were made in \cite{FR} for the values of general $G$-functions by using results of Chudnovsky and Andr\'e and Pad\'e type approximations of the second kind obtained by Siegel's lemma. These studies are motivated by questions on $b$-ary expansions of the function values at $a/b^s$ with integers $s \geq 1$. By applying Theorem 1 with $m = 1$ and denoting $\varphi_1(z) = \varphi(z), \alpha_1 = \alpha = r/s$ and $\alpha_0+\alpha = u/v$, where $(r,s) = (u,v) = 1$, we obtain the following result, where all constants are given explicitly.\\

\noindent
\textbf{Theorem 4}. \textit{Let $a \neq 0$ and $b,B \geq 1$ be integers satisfying
\begin{equation}\label{e5}
b \geq (a_1\left|a\right|)^6, \quad B \leq b^t,
\end{equation}
where $t \geq 0$ and $a_1 = (\vartheta ds_0\epsilon(s_0)\epsilon(v))^{1/4}\tilde{d}\epsilon(s)e^{\theta( s_0/2+s+2v)}$ with any constants $\theta, \vartheta > 1$. Then there exists a positive constant $M_0$ (given explicitly in (\ref{28})) depending on $\alpha_0, \alpha, t, \theta, \vartheta, a$ and $b$ such that, for any $n \in \mathbb{Z}$ and any $M \geq M_0$},
\[
\left|\varphi(\frac{a}{b}) - \frac{n}{Bb^M}\right| \geq \frac{1}{Bb^M(a_1^{18}\left|a\right|^{17})^M}. \\
\]

If $\alpha_0$ is an integer, then $a_1$ can be replaced above by $a_1 = (4\vartheta)^{1/4}e^{3\theta s}$ and in the particular case $\alpha_0 = 1$ by $a_1 = 2^{1/4}e^{3\theta s}$. Concerning the dependence of $M_0$ on $\theta$ and $\vartheta$ above we note that it comes from the use of the prime number theorem 
\[
\pi(x) \leq \theta\frac{x}{\log x}
\]
for all $x \geq c(\theta)$ (as usual $\pi(x)$ denotes the number of primes $\leq x$), and from the use of the inequality
\[
(n + 1)^2 \leq \vartheta^n
\]
for all $n \geq c(\vartheta)$. Further, if $\epsilon, 0 < \epsilon < 1,$ is given and
\[
b^\epsilon > a_1^{18}\left|a\right|^{17},
\]
then the above lower bound implies
\[
\left|\varphi(\frac{a}{b}) - \frac{n}{Bb^M}\right| \geq \frac{1}{Bb^{M(1+\epsilon)}}
\]
for any $n \in \mathbb{Z}$ and any $M \geq \tilde{M}_0$, where $\tilde{M}_0$ depends also on $\epsilon$.

This paper is organized as follows. In the next section we prove a basic lemma, which is then in section 3 used in the explicit construction of Pad\'e approximations of Theorem 1. The denominators of the coefficients of the constructed polynomials are studied and after that the needed upper bounds for the polynomials and remainder terms are obtained in section 4. Theorem 2 follows from more detailed Theorem 5 proved in section 5 and the proof of Theorem 3 is given in the following section. Finally, we shall prove Theorem 4 by using the approximations of Theorem 1 and the main idea of \cite{FR}.

\section{Basic lemma}

Our Pad\'{e} approximation construction is obtained by the following lemma, which is a generalization of \cite[Lemma 1]{F1} considering the case $\alpha_0 = 1$. The result is analogous to \cite[Proposition 2]{So} and we obtain it simply by using the non-vanishing of the determinant $\Delta_N$ below.\\

\noindent
\textbf{Lemma 1.} \textit{Let $\alpha_0$ be a positive number, and assume that $\gamma_1,\ldots,\gamma_m$ are distinct positive numbers. Then the determinant
\[
\Delta_N = \Delta_N(\gamma_1,\ldots,\gamma_N) := \det \left(1 \ \frac{(\gamma_\sigma)_1}{(\gamma_\sigma+\alpha_0)_1} \ \ldots \ \frac{(\gamma_\sigma)_{N-1}}{(\gamma_\sigma+\alpha_0)_{N-1}}\right)_{\sigma=1,\ldots,N} \neq 0,
\]
and the system of linear equations
\begin{equation}\label{1}
1 + \frac{(\gamma_\sigma)_1}{(\gamma_\sigma+\alpha_0)_1}b_1 + \cdots + \frac{(\gamma_\sigma)_N}{(\gamma_\sigma+\alpha_0)_N}b_N = 0, \quad \sigma = 1,\ldots,N,
\end{equation}
has a unique solution
\begin{equation}\label{2}
b_{k+1} = \sum_{\ell=k}^{N-1}(-1)^{\ell+1}\frac{(\alpha_0-1)_{\ell-k}}{(\ell-k)!}\frac{(\alpha_0+\ell+1)_{N-\ell-1}}{(N-\ell-1)!}\prod_{\sigma=1}^N\frac{\gamma_\sigma+\alpha_0+\ell}{\gamma_\sigma}, \quad k = 0,1,\ldots,N-1.
\end{equation}
Moreover, if $\gamma_0 \neq \gamma_\sigma, \sigma = 1,\ldots,N$, is positive and $b_k$ are given in (\ref{2}), then}
\[
1 + \frac{(\gamma_0)_1}{(\gamma_0+\alpha_0)_1}b_1 + \cdots + \frac{(\gamma_0)_N}{(\gamma_0+\alpha_0)_N})b_N \neq 0. \\
\]

\noindent
\textit{Proof.} To prove the non-vanishing of $\Delta_N$ we use induction as in the proof of \cite[Lemma 1]{F1}. Now $\Delta_1(\gamma_1) = 1$. Let $N \geq 2$ and assume that $\Delta_{N-1}(\gamma_1,\ldots,\gamma_{N-1}) \neq 0$ for any distinct positive numbers $\gamma_1,\ldots,\gamma_{N-1}$. We have
\[
\Delta_N(z):= \Delta_N(\gamma_1,\ldots,\gamma_{N-1},z) = \sum_{k=0}^{N-1}\Delta_{Nk}\frac{(z)_k}{(z+\alpha_0)_k} = \frac{d(z)}{(z+\alpha_0)_{N-1}},\]
where $\Delta_{Nk}$ are the cofactors corresponding to the last row of $\Delta_N$ and $d(z)$ is a polynomial of degree at most $N-1$. Since $\Delta_N(\gamma_\sigma) = 0, \ \sigma = 1,\ldots,N-1,$ it follows that
\[
d(z) = c\prod_{\sigma=1}^{N-1}(z - \gamma_\sigma),
\]
with a constant $c$. If $c = 0$, then $\Delta_N(z) \equiv 0$, in particular
\[
\Delta_N(0) = \Delta_{N0} = \pm\frac{\gamma_1\cdots\gamma_{N-1}}{(\gamma_1+\alpha_0)\cdots(\gamma_{N-1}+\alpha_0)}\Delta_{N-1}(\gamma_1+1,\ldots,\gamma_{N-1}+1) = 0
\]
against our induction hypothesis. Therefore $c \neq 0$ and 
\[
\Delta_N = \Delta_N(\gamma_N) = \frac{c}{(\gamma_N+\alpha_0)_{N-1}}\prod_{\sigma=1}^{N-1}(\gamma_N - \gamma_\sigma) \neq 0,
\]
since $\gamma_1,\ldots,\gamma_N$ are distinct.

The coefficient determinant of (\ref{1}) is
\[
\frac{\gamma_1\cdots\gamma_N}{(\gamma_1+\alpha_0)\cdots(\gamma_N+\alpha_0)}\Delta_N(\gamma_1+1,\ldots,\gamma_N+1) \neq 0,
\]
and so (\ref{1}) has a unique solution $b_1,\ldots,b_N$. To determine this solution we define a rational function
\begin{equation}\label{3}
B(x) = 1 + \sum_{k=1}^Nb_k\frac{(x)_k}{(x+\alpha_0)_k} = \frac{b(x)}{(x+\alpha_0)_N},
\end{equation}
where $b(x)$ is a polynomial of degree at most $N$. By (\ref{1}), $b(\gamma_\sigma) = 0, \ \sigma = 1,\ldots,N$, and therefore there exists a constant $c$ such that
\[
b(x) = c\prod_{\sigma=1}^N(x - \gamma_\sigma).
\]
Since $B(0) = 1$, we get
\[
c = (-1)^N(\alpha_0)_N\prod_{\sigma=1}^N\gamma_\sigma^{-1}
\]
and
\begin{equation}\label{4}
B(x) = (-1)^N\frac{(\alpha_0)_N}{(x+\alpha_0)_N}\prod_{\sigma=1}^N\frac{x - \gamma_\sigma}{\gamma_\sigma}.
\end{equation}

For all $k = 0,1,\ldots,N-1$, we may write
\[
\frac{(x)_{k+1}}{(x+\alpha_0)_{k+1}} = 1 + \sum_{\ell=0}^k \frac{c_{k\ell}}{x+\alpha_0+\ell},
\]
where
\[
c_{k\ell} = \lim_{x\rightarrow -\alpha_0-\ell}\frac{(x+\alpha_0+\ell)(x)_{k+1}}{(x+\alpha_0)_{k+1}} = (-1)^{k-\ell+1}\frac{(\alpha_0+\ell-k)_{k+1}}{\ell!(k-\ell)!}, \quad 0\leq\ell\leq k\leq N-1.
\]
Therefore, by (\ref{3}),
\[
B(x) = 1+b_1+\cdots+b_N + \sum_{\ell=0}^{N-1}\sum_{k=\ell}^{N-1}\frac{c_{k\ell}b_{k+1}}{x+\alpha_0+\ell}.
\]
By using (\ref{4}) we now obtain
\[
\sum_{k=\ell}^{N-1}c_{k\ell}b_{k+1} = \lim_{x\rightarrow -\alpha_0-\ell}(x+\alpha_0+\ell)B(x) = (-1)^{\ell}\frac{(\alpha_0)_N}{\ell!(k-\ell)!}\prod_{\sigma=1}^N\frac{\gamma_\sigma+\alpha_0+\ell}{\gamma_\sigma} =: \Gamma_{\ell}
\]
for all $\ell = 0,1,\ldots,N-1$, which gives the system of linear equations
\[
E(b_1,\ldots,b_N)^T = (\Gamma_0/c_{00},\ldots,\Gamma_{N-1}/c_{N-1,N-1})^T,
\] 
where $E = (e_{\ell k})$ is the upper-triangular matrix with
\[
e_{\ell k} = \frac{c_{k\ell}}{c_{\ell\ell}} = (-1)^{k-\ell}\frac{(\alpha_0+\ell-k)_{k-\ell}}{(k-\ell)!}, \quad  0\leq\ell\leq k\leq N-1.
\]
It is proved in \cite[p. 248]{So} that the inverse of $E$ is the upper-triangular matrix $E^{-1} = (d_{\ell k})$, where
\[
d_{\ell k} = \frac{(\alpha_0)_{k-\ell}}{(k-\ell)!}, \quad  0\leq\ell\leq k\leq N-1.
\] 
Therefore
\[
(b_1,\ldots,b_N)^T = E^{-1}(\Gamma_0/c_{00},\ldots,\Gamma_{N-1}/c_{N-1,N-1})^T,
\]
which immediately gives (\ref{2}).

Finally, if also $B(\gamma_0) = 0$, then the system of linear homogeneous equations
\[
b_0 + \frac{(\gamma_\sigma)_1}{(\gamma_\sigma+\alpha_0)_1}b_1 + \cdots + \frac{(\gamma_\sigma)_N}{(\gamma_\sigma+\alpha_0)_N})b_N = 0, \quad \sigma = 0,1,\ldots,N,
\]
would have a solution $b_0 = 1, b_1,\ldots, b_N$, but this is not possible, since the coefficient determinant is $\Delta_{N+1}(\gamma_0,\gamma_1,\ldots,\gamma_N) \neq 0$. Thus necessarily $B(\gamma_0) \neq 0$ and Lemma 1 is proved.

\section{Pad\'e approximations, proof of Theorem 1}

In this section we consider $\varphi_j(z)$ as formal power series in $\mathbb{Q}[[z]]$.
Let $n_1,\ldots, n_m, N_1,\ldots, N_m$ be positive integers and $N = n_1+\cdots+n_m$. We shall construct explicitly type $(N;N_1,\ldots,N_m)$ Pad\'e approximation polynomials of the second kind $Q(z) \neq 0, P_1(z),\ldots, P_m(z)$ for the series $\varphi_1(z),\ldots, \varphi_m(z)$. 
Let
\[
Q(z) = \sum_{k=0}^Na_kz^k.
\]
Then, for each $j=1,\ldots, m$,
\[
Q(z)\varphi_j(z) = \sum_{\mu=0}^\infty c_{j\mu}z^\mu, \quad c_{j\mu} = \sum_{k=0}^{\min\left\{N,\mu\right\}}a_k\frac{(\alpha_j)_{\mu-k}}{(\alpha_j+\alpha_0)_{\mu-k}}, \quad \mu = 0,1,\ldots ,
\]
and
\begin{equation}\label{5}
\ord(Q(z)\varphi_j(z) - P_j(z)) \geq N_j+n_j+1, \quad j = 1,\ldots, m.
\end{equation} 
is satisfied, if
\[
P_j(z) = \sum_{k=0}^{N_j} c_{jk}z^k
\]
and
\begin{equation}\label{6}
c_{j\mu} = 0, \quad \mu = N_j+1,\ldots, N_j+n_j; j = 1,\ldots, m.
\end{equation}
This is a system of $N$ linear homogeneous equations in $N+1$ unknown coefficients $a_k$ and thus $Q(z) \neq 0$ satisfying (\ref{5}) exists.

Assuming $\mu \geq N$ we get
\begin{align*}
&c_{j\mu} = \sum_{k=0}^Na_k\frac{(\alpha_j)_{\mu-k}}{(\alpha_j+\alpha_0)_{\mu-k}} = \\ \frac{(\alpha_j)_{\mu-N}}{(\alpha_j+\alpha_0)_{\mu-N}}&\left(a_N + a_{N-1}\frac{(\alpha_j+\mu-N)_1}{(\alpha_j+\alpha_0+\mu-N)_1} +\ldots+ a_0\frac{(\alpha_j+\mu-N)_N}{(\alpha_j+\alpha_0+\mu-N)_N}\right).
\end{align*}
Therefore, if we denote $a_{N-k} = b_k, \ k=0,1,\ldots, N,$
\begin{align*}
\gamma_1 &= \alpha_1+N_1+1-N, \ldots , \gamma_{n_1} = \alpha_1+N_1+n_1-N, \\
\gamma_{n_1+1} &= \alpha_2+N_2+1-N, \ldots , \gamma_{n_1+n_2} = \alpha_2+N_1+n_2-N, \ldots, \\
\gamma_{n_1+\cdots+n_{m-1}+1} &= \alpha_m+N_m+1-N, \ldots , \gamma_{N} = \alpha_m+N_m+n_m-N,
\end{align*}
and choose $a_N = b_0 = 1, N_j \geq N-1, \ j=1,\ldots,m$, then the system (\ref{6}) is equivalent to the system (\ref{1}) and has a solution (\ref{2}), which has the form
\begin{align}\label{7}
a_N = 1, a_{N-k-1} =& \sum_{\ell=k}^{N-1}(-1)^{\ell+1}\frac{(\alpha_0-1)_{\ell-k}}{(\ell-k)!}\frac{(\alpha_0+\ell+1)_{N-\ell-1}}{(N-\ell-1)!}\prod_{j=1}^m\prod_{\nu=1}^{n_j}\frac{\alpha_j+\alpha_0+N_j+\nu+\ell-N}{\alpha_j+N_j+\nu-N}, \\ &k = 0,1,\ldots,N-1.\notag
\end{align}
Thus Theorem 1 is proved.

We now make $m+1$ constructions by replacing, for each $i=0,1,\ldots, m$ the parameters $N_1,\ldots, N_m$ above by $N_{i1} = N_1+\delta_{i1},\ldots, N_{im} = N_m+\delta_{im}$, in particular $N_{0j} = N_j$ as above. Let the resulting polynomials be, for $i = 0,1,\ldots, m$,
\begin{equation}\label{8}
Q_i(z) = \sum_{k=0}^Na_{ik}z^k, \quad P_{ij}(z) = \sum_{\mu=0}^{N_{ij}} c_{ij\mu}z^\mu, \quad c_{ij\mu} = \sum_{k=0}^{\min\left\{N,\mu\right\}}a_{ik}\frac{(\alpha_j)_{\mu-k}}{(\alpha_j+\alpha_0)_{\mu-k}}, \quad j=1,\ldots, m,
\end{equation} 
and the remainder terms
\[
R_{ij}(z) = Q_i(z)\varphi_j(z) - P_{ij}(z) = \sum_{\mu=N_{ij}+n_j+1}^\infty c_{ij\mu}z^\mu, \quad j=1,\ldots, m.
\]
If $1 \leq i \leq m$, let $\gamma_0 = \alpha_i+N_i+1-N$. Then
\begin{align*}
&c_{ii,N_i+1} = \sum_{k=0}^N a_{ik}\frac{(\alpha_i)_{N_i+1-k}}{(\alpha_i+\alpha_0)_{N_i+1-k}} = \\ \frac{(\alpha_i)_{N_i+1-N}}{(\alpha_i+\alpha_0)_{N_i+1-N}}&\left(a_N + a_{N-1}\frac{(\gamma_0)_1}{(\gamma_0+\alpha_0)_1} +\ldots+ a_0\frac{(\gamma_0)_N}{(\gamma_0+\alpha_0)_N}\right) \neq 0
\end{align*}
by Lemma 1. So $\deg Q_i(z) = N, \deg P_{ij}(z) \leq N_j$ for all $i \neq j$ and $\deg P_{ii}(z) = N_i+1$. The approximations obtained in this way satisfy the following independence lemma.\\

\noindent
\textbf{Lemma 2.} \textit{The determinant
\[
\Omega(z) = \det\left(Q_i(z) \ P_{i1}(z) \ \ldots \ P_{im}(z)\right)_{i=0,1,\ldots,m} = \omega z^{N+N_1+\cdots+N_m+m},
\]
where $\omega \neq 0$ is a constant.}\\

\noindent
\textit{Proof.} The coefficients of the leading terms of $Q_0(z)$ and $P_{ii}(z)$ are 1 and $c_{ii,N_i+1} \neq 0$, respectively. Thus $\Omega(z)$ is a polynomial of exact degree $N+N_1+\cdots+N_m+m$ with the leading coefficient $c_{11,N_1+1}\cdots c_{mm,N_m+1} =: \omega \neq 0$.

On the other hand,
\[
\Omega(z) = (-1)^m\det\left(Q_i(z) \ R_{i1}(z) \ \ldots \ R_{im}(z)\right)_{i=0,1,\ldots,m}.
\]
Here $\ord R_{ij}(z) \geq N_j+n_j+1$, and therefore $\ord \Omega(z) \geq N+N_1+\cdots+N_m+m$, which proves Lemma 2.

\section{Denominators and upper bounds}

The parameters $N_j$ are now chosen in such a way that the obtained approximations can be used to consider the applications studied in Theorems 2 and 4. For this let $n_0 \geq \max \left\{n_1,\ldots,n_m\right\}$ and fix $N_j = N+n_0-n_j =: \tilde{N}-n_j, \ j=1,\ldots, m.$ By substituting then $N_{ij}$ to (\ref{e3}) we get, for each $i = 0,1,\ldots, m$, 
\begin{align}\label{9}
&a_{iN} = 1, a_{i,N-k-1} = \sum_{\ell=k}^{N-1}(-1)^{\ell+1}\frac{(\alpha_0-1)_{\ell-k}}{(\ell-k)!}\frac{(\alpha_0+\ell+1)_{N-\ell-1}}{(N-\ell-1)!}\times \\ &\prod_{j=1}^m\frac{(\alpha_j+\alpha_0+n_0-n_j+\delta_{ij}+\ell+1)_{n_j}}{(\alpha_j+n_0-n_j+\delta_{ij}+1)_{n_j}}, \ k = 0,1,\ldots,N-1.\notag
\end{align}
To consider these coefficients we first give a lemma from \cite [pp. 145-147]{M} considering the quotients
\[
\frac{(\alpha+1)_n}{n!} =: \frac{p_n}{q_n}, \ (p_n,q_n) = 1, q_n \geq 1, n = 0,1,\ldots,
\]
where $\alpha = r/s \neq -1,-2,\ldots$ with integers $r$ and $s \geq 1, (r,s) = 1$.\\

\noindent
\textbf{Lemma 3}. \textit{Let
\[
U_n = \prod_{p\nmid s}p^{M_p}, \ V_n = s\prod_{p\mid s}p^{\nu_p(n)},
\]
where 
\[
M_p = [\log \max \{\left|r+s\right|, \left|r+sn\right|\}/\log p] \leq [\log(\left|r\right|+sn)/\log p], \quad \nu_p(n) = \sum_{t=1}^\infty\left[\frac{n}{p^t}\right] \leq \frac{n}{p-1}.
\]
Then the least common multiples of $p_0,p_1,\ldots,p_n$ and of $q_0,q_1,\ldots,q_n$ are divisors of $U_n$ and $V_n$, respectively.}\\

To consider the above $a_{ik}$ we first note that
\begin{align*}
\frac{(\alpha_j+\alpha_0+n_0-n_j+\delta_{ij}+\ell+1)_{n_j}}{(\alpha_j+n_0-n_j+\delta_{ij}+1)_{n_j}} &= \frac{(u_j+(n_0-n_j+\delta_{ij}+\ell+1)v_j)\cdots(u_j+(n_0+\delta_{ij}+\ell)v_j)}{n_j!s_0^{n_j}} \\ & \times\frac{n_j!d_j^{n_j}}{(r_j+(n_0-n_j+\delta_{ij}+1)s_j)\cdots(r_j+(n_0+\delta_{ij})s_j}.  
\end{align*} 
Therefore Lemma 3 immediately implies that
\[
D_1(N)a_{ik} \in \mathbb{Z}, \quad i=0,1,\ldots,m; k = 0,1,\ldots,N,
\]
if
\[
D_1(N) = s_0^{2N-1}\prod_{p\mid s_0}p^{\nu_p(N-1)}\times\prod_{j=1}^m\left\{\prod_{p\mid v_j}p^{\nu_p(n_j)}\prod_{p\nmid s_j}p^{M_p(r_j+(n_0+1)s_j)}\right\},
\]
where we denote $M_p(x) = \left[\log x/\log p\right]$ for all $x > 0$. So we have
\begin{equation}\label{10}
D_1(N) \leq s_0^{2N-1}\epsilon(s_0)^{N-1}\prod_p p^{mM_p(R+S+Sn_0)}\times\prod_{j=1}^m\epsilon(v_j)^{n_j} \leq (s_0^2\epsilon(s_0)\epsilon(v))^N\prod_p p^{mM_p(R+S+Sn_0)}.
\end{equation}

For the consideration of the coefficients of $P_{ij}(z)$ we see that
\begin{equation}\label{11}
\frac{(\alpha_j)_{\mu-k}}{(\alpha_j+\alpha_0)_{\mu-k}} = \frac{r_j(r_j+s_j)\cdots(r_j+s_j(\mu-k-1))}{(\mu-k)!d_j^{\mu-k}}\frac{(\mu-k)!s_0^{\mu-k}}{u_j(u_j+v_j)\cdots(u_j+v_j(\mu-k-1))},
\end{equation}
and therefore, by Lemma 3,
\begin{equation}\label{e11}
\left(\frac{d_j}{(d_j,s_0)}\right)^\mu \prod_{p\mid s_j}p^{\nu_p(\mu)}\times\prod_{p\nmid v_j}p^{M_p(u_j+v_j\mu)}\times\frac{(\alpha_j)_{\mu-k}}{(\alpha_j+\alpha_0)_{\mu-k}} \in \mathbb{Z}
\end{equation}
for all $k = 0,1,\ldots,\mu; \mu \geq 0$. Thus we get the following lemma, where
\begin{equation}\label{12}
D_2(N) = \left(\frac{d}{(d,s_0)}\right)^{\tilde{N}}\prod_{p\mid s}p^{\nu_p(\tilde{N})}\times\prod_pp^{M_p(U+V\tilde{N})} \leq \left(\tilde{d}\epsilon(s)\right)^{\tilde{N}}\prod_pp^{M_p(U+V\tilde{N})}.\\
\end{equation}

\noindent
\textbf{Lemma 4}. \textit{If $D(N) = D_1(N)D_2(N)$, then}
\[
D(N)Q_i(z), D(N)P_{ij}(z) \in \mathbb{Z}\left[z\right], \quad i = 0,1,\ldots,m; j = 1,\ldots,m.\\
\]

By the prime number theorem, for any $\theta > 1$, the number of primes $p \leq x$ satisfies
\[
\pi(x) \leq \theta\frac{x}{\log x}
\]
for all $x\geq c(\theta)$, in particular we may take $c(8\log 2) = 2$, see for example \cite[p. 296]{Bu}. Therefore (\ref{10}) and (\ref{12}) give, for all $\min\left\{n_0,N\right\} \geq c(\theta)$,
\begin{align}\label{13}
&D(N) \leq e^{c_1+c_2n_0+c_3N+c_4\tilde{N}},\\ c_1 = \theta(m(R+S)+U), \ c_2 &= m\theta S, \ c_3 = \log(s_0^2\epsilon(s_0)\epsilon(v)), \ c_4 = \theta V+\log(d\epsilon(s)/(d,s_0)) \notag.
\end{align}

In the same way by using Lemma 3 we also have
\begin{align*}
\left|a_{ik}\right| &\leq N\prod_{p\nmid s_0}p^{2M_p(r_0+s_0(N-1))}\times\prod_{j=1}^m\left\{\frac{d_j^{n_j}}{s_0^{n_j}}\prod_{p\mid s_j}p^{\nu_p(n_j)}\times\prod_{p\nmid v_j}p^{M_p(U+V\tilde{N})}\right\} \\ &\leq Ns_0^{-N}\prod_pp^{2M_p(r_0+s_0(N-1))+mM_p(U+V\tilde{N})}\times\prod_{j=1}^m(d_j\epsilon(s_j))^{n_j}, 
\end{align*}
which gives
\begin{equation}\label{14}
\left|a_{ik}\right| \leq Ne^{c_5+c_6N+c_7\tilde{N}}, \quad c_5 = \theta(2(r_0-s_0)+mU), \ c_6 = 2\theta s_0+\log(d\epsilon(s)/s_0), \ c_7 = m\theta V
\end{equation}
for all $i=0,1,\ldots,m; k=0,1,\ldots,N$ and $N \geq c(\theta)$. Thus the following lemma holds.\\

\noindent
\textbf{Lemma 5}. \textit{For all $N \geq c(\theta)$ and $\left|z\right| \geq 2$},
\begin{align*}
\left|Q_i(z)\right| \leq 2Ne^{c_5+c_6N+c_7\tilde{N}}&\left|z\right|^N, \ \left|P_{ij}(z)\right| \leq 2N(N+1)e^{c_5+c_6N+c_7\tilde{N}}\left|z\right|^{\tilde{N}-n_j+1},\\
&i=0,1,\ldots,m; \ j=1,\ldots,m.\\
\end{align*}

For the proof of Theorem 2 we still need to estimate the $p$-adic values of the remainder terms $R_{ij}(\beta)$, where a rational $\beta = a/b \neq 0$ satisfies $\left|a\right|_p < 2^{-\delta(2,p)}\left|s\right|_p$. We note first that, for all $\mu \geq \tilde{N}$,
\[
\left|D_1(N)c_{ij\mu}\right|_p \leq \max_{0\leq k \leq N}\left|\frac{D_1(N)a_{ik}(\alpha_j)_{\mu-k}}{(\alpha_j+\alpha_0)_{\mu-k}}\right|_p \leq \max_{0\leq k \leq N}\left|\frac{(\alpha_j)_{\mu-k}}{(\alpha_j+\alpha_0)_{\mu-k}}\right|_p,
\]
since $D_1(N)a_{ik} \in \mathbb{Z}$. Therefore, by (\ref{11}), (\ref{e11}) and Lemma 3,
\[
\left|D_1(N)c_{ij\mu}\right|_p \leq \left|\tilde{d}\right|_p^{-\mu}p^{\delta(p)\nu_p(\mu)+M_p(U+V\mu)} \leq \left|\tilde{d}\right|_p^{-\mu}p^{\delta(p)\nu_p(\mu)}(U+V\mu),
\]
where $\delta(p) = 1$, if $p\mid s$, and $\delta(p) = 0$ otherwise. This gives
\begin{align*}
\left|D_1(N)R_{ij}(\beta)\right|_p &= \left|\sum_{\mu=\tilde{N}+1}^\infty D_1(N)c_{ij\mu}a^\mu \right|_p \leq \max_{\mu\geq\tilde{N}+1}\left\{(U+V\mu)p^{\delta(p)\mu/(p-1)}\left|\frac{a}{\tilde{d}}\right|_p^\mu\right\} \\ &= (U+V(\tilde{N}+1))p^{\delta(p)(\tilde{N}+1)/(p-1)}\left|\frac{a}{\tilde{d}}\right|_p^{\tilde{N}+1},
\end{align*}
since $\left|a/\tilde{d}\right|_pp^{\delta(p)/(p-1)} < 1/\sqrt{p}$ by the above assumption $\left|a\right|_p < 2^{-\delta(2,p)}\left|s\right|_p$ and the fact $\tilde{d} \mid s$. Thus, by (\ref{12}),
\begin{align*}
\left|D(N)R_{ij}(\beta)\right|_p &\leq \left|D_2(N)\right|_p(U+V(\tilde{N}+1))p^{\delta(p)(\tilde{N}+1)/(p-1)}\left|\frac{a}{\tilde{d}}\right|_p^{\tilde{N}+1} \\ &\leq \tilde{d}\frac{U+V(\tilde{N}+1)}{U+V\tilde{N}}p^{\delta(p)((\tilde{N}+1)/(p-1)-\nu_p(\tilde{N}))}\left|a\right|_p^{\tilde{N}+1}.
\end{align*}
Here
\[
\frac{\tilde{N}+1}{p-1} - \nu_p(\tilde{N}) \leq t+4,
\]
where $t$ is the integer satisfying $p^t \leq \tilde{N} < p^{t+1}$. This implies
\begin{equation}\label{a14}
\left|D(N)R_{ij}(\beta)\right|_p \leq 2\tilde{d}\left|a\right|^{4\delta(p)}\tilde{N}^{\delta(p)}\left|a\right|_p^{\tilde{N}+1}.
\end{equation}
By combining this estimate with Lemmas 2, 4 and 5 we obtain the following result considering the integers $Q_i, P_{ij}$ and the $p$-adic numbers $R_{ij}$ defined by
\[
Q_i := D(N)b^{\tilde{N}}Q_i(\beta), \ P_{ij} := D(N)b^{\tilde{N}}P_{ij}(\beta), \ R_{ij} := D(N)b^{\tilde{N}}R_{ij}(\beta), \ i=0,1,\ldots,m; j=1,\ldots,m.\\
\]

\noindent
\textbf{Lemma 6}. \textit{Let a rational $\beta = a/b \neq 0$ satisfy $\left|\beta\right| \geq 2$ and $\left|a\right|_p < 2^{-\delta(2,p)}\left|s\right|_p$. Then the above defined $Q_i, P_{ij} \in \mathbb{Z}, \ i=0,1,\ldots,m; j=1,\ldots,m,$ and satisfy
\[
\det (Q_i \ P_{i1} \ \ldots, \ P_{im})_{i=0,1,\ldots,m} \neq 0.
\]
Further, for all $\tilde{N} \geq \tilde{N}_1 := \max \{(m+1)c(\theta), \ c_1+c_5, \ \log(2\tilde{d})+4\delta(p)\log\left|a\right|\}$, we have 
\[
\left|Q_i\right| \leq e^{c_2n_0+c_8\tilde{N}}b^{n_0}\left|a\right|^{\tilde{N}-n_0}, \  \left|P_{ij}\right| \leq e^{c_2n_0+c_8\tilde{N}}b^{n_j}\left|a\right|^{\tilde{N}-n_j+1}, 
\]
\[
\left|R_{ij}\right|_p \leq e^{2\tilde{N}}\left|a\right|_p^{\tilde{N}+1}, \
i=0,1,\ldots,m; j=1,\ldots,m,
\]
where $c_8 = c_3+c_4+c_6+c_7+3$.}\\

\section{Proof of Theorem 2}

We now consider a linear form $\ell_0+\ell_1\varphi_1(\beta)+\cdots+\ell_m\varphi_m(\beta)$ in $p$-adic numbers $\varphi_j(\beta)$, where $\beta$ is given in Lemma 6 and $(\ell_0,\ell_1,\ldots,\ell_m)\in \mathbb{Z}^{m+1}\setminus\{\underline{0}\}$, and denote, as in Theorem 2, $h_j=\max \{1,\left|\ell_j\right|\} \ (j=1,\ldots,m), h_0=\max \{\left|\ell_0\right|,\left|\ell_1\right|,\ldots,\left|\ell_m\right|\}$ and $\tilde{H}=\prod_{j=0}^mh_j$.\\

\noindent
\textbf{Theorem 5}. \textit{Let $\alpha_0,\alpha_1,\ldots, \alpha_m$ satisfy the assumptions of Theorem 2, and assume that $\tau > 0, \delta \geq 0$ and $\beta$ satisfy
\begin{equation}\label{15}
\tau > 4\delta(1+(m+1)\tau).
\end{equation}
and  
\begin{equation}\label{16}
(a,b) = 1, \left|a\right|_p \leq \min\left\{2^{-\delta(2,p)}\left|s\right|_p, \left|a\right|^{-1+\delta}\right\}, \ \left|a\right| > b^{2(1+1/\tau)}e^{2(c_2(1+1/\tau)+(c_8+2)(m+1+1/\tau))}.
\end{equation}
Then
\[
\left|\ell_0+\ell_1\varphi_1(\beta)+\cdots+\ell_m\varphi_m(\beta)\right|_p > \tilde{H}^{-1-(m+1)\tau}
\]
for all $(\ell_0,\ell_1,\ldots,\ell_m) \in \mathbb{Z}^{m+1}\setminus\{\underline{0}\}$ and $\tilde{H} \geq \tilde{H}_0$, where} 
\begin{equation}\label{17}
\log \tilde{H}_0 = \max \left\{\frac{(\tilde{N}_1+m+1)\log \left|a\right|}{1+(m+1)\tau}, \frac{8\log \left|a\right|}{\tau}\right\}.\\
\end{equation}

\noindent
\textit{Proof}. Let $L=\ell_0+\ell_1\varphi_1(\beta)+\cdots+\ell_m\varphi_m(\beta)$. By Lemma 6 there exists an index $i, 0\leq i\leq m$, such that the integer $\Lambda=Q_i\ell_0+P_{i1}\ell_1+\cdots+P_{im}\ell_m \neq 0$. Now
\begin{equation}\label{18}
Q_iL = \Lambda + \sum_{j=1}^m R_{ij}\ell_j.
\end{equation} 
We next choose
\[
n_j = \left[\frac{\log(h_j\tilde{H}^\tau)}{\log \left|a\right|}\right], \quad j=0,1,\ldots,m,
\]
and prove that then
\begin{equation}\label{19}
\left|\Lambda\right|_p > \left|\sum_{j=1}^m R_{ij}\ell_j\right|_p,
\end{equation}
if (\ref{17}) holds.

By Lemma 6 and our choice of $n_j$,
\[
\left|\Lambda\right| \leq e^{c_2n_0+c_8\tilde{N}}\left|a\right|^{\tilde{N}+1}\sum_{j=0}^m b^{n_j}h_j\left|a\right|^{-n_j} < e^{c_2n_0+(c_8+1)\tilde{N}}\left|a\right|^{\tilde{N}+2}\tilde{H}^{-\tau+\frac{(1+\tau)\log b}{\log \left|a\right|}},
\]
and so
\begin{equation}\label{20}
\left|\Lambda\right|_p > e^{-c_2n_0-(c_8+1)\tilde{N}}\left|a\right|^{-\tilde{N}-2}\tilde{H}^{\tau-\frac{(1+\tau)\log b}{\log \left|a\right|}}.
\end{equation}
On the other hand, by Lemma 6 and (\ref{16}),
\[
\left|\sum_{j=1}^m R_{ij}\ell_j\right|_p \leq \max \{\left|R_{ij}\right|_p\} \leq e^{2\tilde{N}}\left|a\right|_p^{\tilde{N}+1} \leq e^{2\tilde{N}}\left|a\right|^{(-1+\delta)(\tilde{N}+1)}.
\]
Therefore (\ref{19}) follows if
\[
\tilde{H}^\tau > \tilde{H}^{\frac{(1+\tau)\log b}{\log \left|a\right|}}e^{c_2n_0+(c_8+2)\tilde{N}}\left|a\right|^{1+\delta(\tilde{N}+1)}.
\]
By our choice of $n_j$,
\[
\tilde{N} = n_0+n_1+\cdots+n_m \leq \frac{(1+(m+1)\tau)\log \tilde{H}}{\log \left|a\right|}, \ n_0 \leq \frac{(1+\tau)\log \tilde{H}}{\log \left|a\right|},
\]
and so the above inequality follows from
\[
\tau > \frac{(c_2+\log b)(1+\tau)+(c_8+2)(1+(m+1)\tau)}{\log \left|a\right|}+\delta(1+(m+1)\tau)+\frac{2\log \left|a\right|}{\log \tilde{H}},
\]
which is a consequence of (\ref{15}), (\ref{16}) and (\ref{17}). Note that for the condition $\tilde{N} \geq \tilde{N}_1$ of Lemma 6 we need to assume here that
\[
\log \tilde{H} \geq \frac{(\tilde{N}_1+m+1)\log \left|a\right|}{1+(m+1)\tau}.
\]

By using (\ref{18}), (\ref{19}) and (\ref{20}) we now get
\[
\left|L\right|_p \geq \left|Q_iL\right|_p = \left|\Lambda\right|_p >e^{-c_2n_0-(c_8+1)\tilde{N}}\left|a\right|^{-\tilde{N}-2}\tilde{H}^{\tau-\frac{(1+\tau)\log b}{\log \left|a\right|}} \geq \tilde{H}^{-\eta},
\]
where, by (\ref{16}) and (\ref{17}),
\[
\eta = 1+(m+1)\tau - \tau + \frac{(1+\tau)(c_2+\log b)}{\log \left|a\right|} + \frac{(c_8+1)(1+(m+1)\tau)}{\log \left|a\right|} + \frac{2\log \left|a\right|}{\log \tilde{H}} \leq 1+(m+1)\tau.
\]
This proves Theorem 5.

The choice $\tau = \epsilon/(m+1)$ and $\delta = \epsilon/(8(m+1))$ in Theorem 5 implies Theorem 2.

\section{Proof of Theorem 3}

Here we use Lemma 6 and (\ref{a14}), where we now assume that $n_0=n_1=\cdots=n_m=:n$ and $b=1, (a,s) = 1$. Then, for $p \mid a$ and all $n \geq \tilde{N}_1/(m+1)$,
\[
\left|Q_i\right|, \left|P_{ij}\right| \leq e^{c_9n}\left|a\right|^{mn+1}, \ \left|R_{ij}\right|_p \leq 2\tilde{d}\left|a\right|_p^{(m+1)n+1},
\]
where $c_9 = c_2+(m+1)c_8$.

Suppose now that an $a$-global relation
\[
L = \ell_0 + \ell_1\varphi_1(a) + \cdots + \ell_m\varphi_m(a) = 0
\]
holds and denote $h= \max \{\left|\ell_0\right|,\left|\ell_1\right|,\ldots,\left|\ell_m\right|\}$. As in the proof of Theorem 2 there exists an integer $i, 1\leq i\leq m,$ such that the integer $\Lambda = Q_i\ell_0+P_{i1}\ell_1+\cdots+P_{im}\ell_m \neq 0$. Since $L = 0$ for all $p \mid a$, we have
\[
\Lambda = -\sum_{j=1}^m R_{ij}\ell_j
\]
for all $p \mid a$. Further, the number of these primes is $\leq \log \left|a\right|/\log 2$. Thus
\begin{align*}
&1 = \left|\Lambda\right|\prod_{p}\left|\Lambda\right|_p \leq \left|\Lambda\right|\prod_{p\mid a}\left|\Lambda\right|_p = \left|\Lambda\right|\prod_{p\mid a}\left|\sum_{j=1}^m R_{ij}\ell_j\right|_p \leq (m+1)he^{c_9n}\left|a\right|^{mn+1}\times \\
&\prod_{p\mid a}2\tilde{d}\left|a\right|_p^{(m+1)n+1} \leq (m+1)h(2\tilde{d})^{\log \left|a\right|/\log 2}e^{(c_9-\log \left|a\right|)n}.
\end{align*}
If $\log \left|a\right|>c_9$, then we get a contradiction for all sufficiently large $n$. Therefore $a$-global relations are possible only if $\left|a\right| \leq C$, where $\log C = c_9$. This proves Theorem 3, since 
\[
c_9 = m\theta S + (m+1)(3 + \log(d\tilde{d}s_0\epsilon(s_0)\epsilon(s)^2\epsilon(v)) + \theta(2s_0 + (m+1)V))
\]
and we may choose $\theta$ arbitrarily close to 1.\\

\section{Restricted approximations, proof of Theorem 4}

In the consideration of restricted approximations we suppose the notations and assumptions to be as in Theorem 4 and apply the main idea of \cite{FR}. Since $m = 1$, we have $N = n_1$ and the coefficients of $Q_i(z), \ i = 0,1,$ given in (\ref{9}) are now
\begin{align*}
a_{in_1} = 1, a_{i,n_1-k-1} = \sum_{\ell=k}^{n_1-1}(-1)^{\ell+1}&\frac{(\alpha_0-1)_{\ell-k}}{(\ell-k)!}\frac{(\alpha_0+\ell+1)_{n_1-\ell-1}}{(n_1-\ell-1)!}\frac{(\alpha+\alpha_0+n_0-n_1+\delta_{ij}+\ell+1)_{n_1}}{(\alpha+n_0-n_1+\delta_{ij}+1)_{n_1}}, \\  &k = 0,1,\ldots,n_1-1.
\end{align*}
Similarly to the above studies leading to (\ref{10}) and (\ref{12}) we have
\[
D_1Q_i(z), \quad D_1D_2P_i(z) \in \mathbb{Z}[z], \ i=0,1,
\]
where, for all $n_0 \geq c(\theta)$,
\begin{equation}\label{21}
D_1 = s_0^{2n_1-1}\prod_{p\mid s_0}p^{\nu_p(n_1-1)}\times\prod_{p\mid v_j}p^{\nu_p(n_1)}\times\prod_{p\nmid s}p^{M_p(r+s(n_0+1))} \leq (s_0^2\epsilon(s_0)\epsilon(v))^{n_1}e^{\theta(r+s(n_0+1))},
\end{equation}
and
\begin{equation}\label{22}
D_2 = \tilde{d}^{n_0+1}\prod_{p\mid s}p^{\nu_p(n_0+1)}\times\prod_{p\nmid v}p^{M_p(u+vn_0)} \leq (\tilde{d}\epsilon(s))^{n_0+1}e^{\theta(u+vn_0)}.
\end{equation}
By the construction $\deg Q_i(z) = n_1, \deg P_0(z) \leq n_0, \deg P_1(z) = n_0+1$.

We also have, or all $n_1 \geq c(\theta)$,
\begin{align*}
\left|a_{ik}\right| &\leq n_1d^{n_1}s_0^{-n_1}\prod_{p\nmid s_0} p^{2M_p(r_0+s_0n_1)}\times\prod_{p\mid s}p^{\nu_p(n_1)}\times\prod_{p\nmid v} p^{M_p(u+v\tilde{N})} \\ &\leq n_1e^{\theta(2r_0+u)}(ds_0^{-1}\epsilon(s))^{n_1}e^{\theta(2s_0n_1+v\tilde{N})} =: E_1,
\end{align*}
where now $\tilde{N} = n_0+n_1$. Thus we have, for all $\left|z\right| < 1, i = 0,1$,
\begin{equation}\label{23}
\left|Q_i(z)\right| \leq E_1\frac{1}{1-\left|z\right|}, \ \left|R_i(z)\right| \leq (n_1+1)E_1\left|z\right|^{\tilde{N}+1}\frac{1}{1-\left|z\right|},
\end{equation}
if $n_1 \geq c(\theta)$.

Let $a, b, n, B$ and $M$ be as in Theorem 4. By Lemma 2 there exists an $i \in \{0,1\}$ such that
\[
nQ_i(\frac{a}{b}) - Bb^MP_i(\frac{a}{b}) \neq 0.
\]
Here $D_1b^{n_1}Q_i(a/b) =: U_i \in \mathbb{Z}$ and $D_1D_2b^{n_0+1}P_i(a/b) =: V_i \in \mathbb{Z}$, and therefore
\[
0 \neq D_1D_2b^{n_0+1}\left(nQ_i(\frac{a}{b}) - Bb^MP_i(\frac{a}{b})\right) = D_2b^{n_0-n_1+1}U_i - Bb^MV_i =: W_i \in \mathbb{Z}.
\]
If $n_0-n_1+1 \geq M$, then $b^M \mid W_i$ and
\[
b^M \leq \left|W_i\right| \leq D_1D_2b^{n_0+1}\left(\left|Q_i(\frac{a}{b})\right|\left|n - Bb^M\varphi(\frac{a}{b})\right| + Bb^M\left|Q_i(\frac{a}{b})\varphi(\frac{a}{b}) - P_i(\frac{a}{b})\right|\right).
\]
implying
\[
\left|Q_i(\frac{a}{b})\right|\left|n - Bb^M\varphi(\frac{a}{b})\right| \geq \frac{b^M}{D_1D_2b^{n_0+1}} - Bb^M\left|R_i(\frac{a}{b})\right|.
\]
This inequality gives a lower bound
\begin{equation}\label{24}
\left|Q_i(\frac{a}{b})\right|\left|n - Bb^M\varphi(\frac{a}{b})\right| \geq \frac{b^M}{2D_1D_2b^{n_0+1}},
\end{equation}
if
\[
B\left|R_i(\frac{a}{b})\right| \leq \frac{1}{2D_1D_2b^{n_0+1}}.
\]
By (\ref{21}), (\ref{22}) and (\ref{23}) this inequality holds for all $n_0 \geq c(\theta)$, if
\[
4n_1(n_1+1)\tilde{d}\epsilon(s)\left|a\right|e^{\theta(r+s+2r_0+2u)}(ds_0\epsilon(s_0)\epsilon(v))^{n_1}\tilde{d}^{n_0}\epsilon(s)^{\tilde{N}}e^{\theta(2s_0n_1+(s+v)n_0+v\tilde{N})}\left|a\right|^{\tilde{N}}Bb^{-n_1} \leq 1.
\]
If $\vartheta>1$, then there exists a positive constant $c(\vartheta)$ such that $(n_1+1)^2 < \vartheta^{n_1}$ for all $n_1 \geq c(\vartheta)$. So the above equality is satisfied, if $n_1 \geq \max\left\{c(\theta), c(\vartheta)\right\}$ and
\begin{equation}\label{25}
a_2\left|a\right|(\vartheta ds_0\epsilon(s_0)\epsilon(v))^{n_1}\tilde{d}^{n_0}\epsilon(s)^{\tilde{N}}e^{\theta(2s_0n_1+(s+v)n_0+v\tilde{N})}\left|a\right|^{\tilde{N}}Bb^{-n_1} \leq 1,
\end{equation}
where we denote
\[
a_2 = 4\tilde{d}\epsilon(s)e^{\theta(r+s+2r_0+2u)}.
\] 

We now choose $n_1 = h, n_0 = \left[xh\right]$, where
\[
x = \frac{\log b}{2\log(a_1\left|a\right|)}.
\]
Then $x \geq 3$ by (\ref{e5}), and therefore (\ref{25}) follows from
\[
a_2\left|a\right|(a_1\left|a\right|)^{(x+1)h}Bb^{-h} \leq 1.
\]
Here
\begin{equation}\label{26}
(a_1\left|a\right|)^{x+1} = b^{1/2}e^{\log(a_1\left|a\right|)} \leq b^{2/3},
\end{equation}
where we again used (\ref{e5}). If $B \leq b^t$ and $h \geq 4t$, it is thus enough to have 
\[
a_2\left|a\right|b^{-h/12} \leq 1
\]
or
\[
h \geq 12\frac{\log(a_2\left|a\right|)}{\log b}
\]
for the validity of (\ref{25}), which implies (\ref{24}). Since we need $n_0-n_1+1 \geq M$ for (\ref{24}), the above means that (\ref{24}) holds if $h$ satisfies
\begin{equation}\label{27}
h \geq \max\left\{12\frac{\log(a_2\left|a\right|)}{\log b}, \ 4t, \ \frac{M}{x-1}, \ c(\theta), \ c(\vartheta), \ 4\right\}.
\end{equation}                                                                                                      

We now fix the parameter h by
\[
h = \left[\frac{M}{x-2}\right].
\] 
For the condition $h \geq M/(x-1)$ we should have
\[
\frac{M}{x-2} - \frac{M}{x-1} \geq 1.
\]
Thus (\ref{27}) follows, if $M \geq M_0$, where
\begin{equation}\label{28}
M_0 := \frac{\log b}{\log(a_1\left|a\right|)}\max\left\{6\frac{\log(a_2\left|a\right|)}{\log b}+\frac{1}{2}, \ \frac{4t+1}{2}, \ \frac{1}{4}\left(\frac{\log b}{\log(a_1\left|a\right|)}\right), \ \frac{1+\max\left\{c(\theta),c(\vartheta),4\right\}}{2}\right\}.
\end{equation} 

Now, by (\ref{21}), (\ref{22}), (\ref{23}) and (\ref{26}),
\[
2D_1D_2b^{n_0+1}\left|Q_i(\frac{a}{b})\right| \leq a_2b\left(\frac{(a_1\left|a\right|)^{1+x}b^x}{\left|a\right|^{1+x}}\right)^h \leq a_2b^{1-h/3}(\frac{b}{\left|a\right|})^{(1+x)h} \leq (\frac{b}{\left|a\right|})^{\frac{x+1}{x-2}M}.
\]
Further, since $x \geq 3$,
\[
\frac{x+1}{x-2} = 1 + \frac{3}{x-2} \leq 1 + \frac{9}{x} = 1 + 18\frac{\log(a_1\left|a\right|)}{\log b},
\] 
and so
\[
(\frac{b}{\left|a\right|})^{\frac{x+1}{x-2}M} \leq (a_1^{18}\left|a\right|^{17})^Mb^M.
\]
Hence, by using (\ref{24}), we get
\[
\left|\varphi(\frac{a}{b}) - \frac{n}{Bb^M}\right| \geq \frac{1}{Bb^M(a_1^{18}\left|a\right|^{17})^M},
\]
which proves Theorem 4.

\begin{small}

\noindent
Keijo V\"a\"an\"anen \\
Department of Mathematical Sciences  \\
University of Oulu \\
P. O. Box 3000 \\
90014 Oulun yliopisto, Finland \\
E-mail: keijo.vaananen@oulu.fi \\

\end{small}


\begin{center}
\textbf{References}
\end{center} 
\begin{thebibliography}{}
\bibitem{B} A. Baker, \textit{On some Diophantine inequalities involving the exponential function}, Canad. J. Math. 17 (1965), 616-626.

\bibitem{Bo} E. Bombieri, \textit{On G-functions}, Recent Progress in Analytic Number Theory 2, Acad. Press (1981), 1-67.

\bibitem{Bu} P. Bundschuh, \textit{Einf\"uhrung in die Zahlentheorie}, 4 Aufl., Springer-Lehrbuch, Springer, 1998.

\bibitem{ELM} A.-M. Ernvall-Hyt\"onen, K. Lepp\"al\"a, T. Matala-aho, \textit{An explicit Baker-type lower bound of exponential values}, Proc. Roy. Soc. Edinburgh Sect. A 145 (2015), 1153-1181.

\bibitem{F} N. I. Fel'dman, \textit{Lower estimates for some linear forms}, Vestnik Moscov. Univ. Ser. I, Mat. Meh. 22, No. 2 (1967), 63-72.

\bibitem{F1} N. I. Fel'dman, \textit{On a linear form}, Acta Arith. 21 (1972), 347-355.

\bibitem{FR} S. Fischler and T. Rivoal, \textit{Rational approximation to values of G-functions, and their expansions in integer bases}, Manuscripta math. 155, 3-4 (2018), 579-595.  

\bibitem{L} L. Leinonen, \textit{A Baker-type linear independence measure for the values of generalized Heine series}, J. Algebra Number Theory Acad. 4 (2014), 49-75.

\bibitem{M1} K. Mahler, \textit{On a paper by A. Baker on the approximation of rational powers of e}, Acta Arith. 27 (1975), 61-87.

\bibitem{M2} K. Mahler, \textit{Lectures on Transcendental Numbers}, Lecture Notes in Mathematics 546, Springer, 1976.

\bibitem{M} Ju. N. Makarov, \textit{On the estimate of the measure of linear independence for the values of E-functions}, Vestnik Moscov. Univ. Ser. I, Mat. Meh. 33, No. 2 (1978), 3-12.

\bibitem{S} O. Sankilampi, \textit{On the linear independence measures of the values of some q-hypergeometric and hypergeometric functions and some applications}, PhD thesis, Univ.of Oulu, 2006.

\bibitem{So} V. N. Sorokin, \textit{On the irrationality of the values of hypergeometric functions}, Sb. Math. 55 (1986), 243-257.

\bibitem{V} K. V\"a\"an\"anen, \textit{On linear forms of a certain class of G-functions and p-adic G-functions}, Acta Arith. 36 (1980), 273-295.

\bibitem{V1} K. V\"a\"an\"anen, \textit{On a result of Fel'dman on linear forms in the values of some $E$-functions}, Ramanujan J., to appear.

\bibitem{VZ} K. V\"a\"an\"anen and W. Zudilin, \textit{Baker-type estimates for linear forms in the values of q-series}, Canad. Math. Bull. 48 (2005), 147-160.

\bibitem{Z} W. Zudilin, \textit{Lower bounds for polynomials in the values of certain entire functions}, Sb. Math. 187 (1996), 1791-1818.

\end{thebibliography}
\end{document}